\documentclass[11pt]{article}
\usepackage{graphicx} 
\usepackage{url}
\usepackage{enumitem}
\usepackage{amsmath}
\usepackage{amssymb}
\usepackage{amsfonts}
\usepackage{breakcites}
\usepackage[fleqn]{nccmath}
\usepackage{ragged2e}
\usepackage{cite}
\usepackage{hyperref}

\usepackage{mathtools}
\newtheorem{Thm}{Theorem}

\newtheorem{Prop}{Proposition}

\newtheorem{Rem}{Remark}

\renewcommand{\bf}[1]{\textbf{#1}}

\newcommand{\equa}{\begin{equation}}
\newcommand{\equb}{\end{equation}}

\newcommand{\eqa}{\begin{equation*}}
\newcommand{\eqb}{\begin{equation*}}
\newcommand{\eqn}{\begin{eqnarray*}}
\newcommand{\eqnb}{\end{eqnarray*}}

\usepackage{lmodern} 
\usepackage{graphicx,scalerel}

\usepackage{indentfirst}

\begin{document}
\begin{center}
\LARGE{\bf{Chaotic and Predictable Representations for Markov Additive Processes with Lévy Modulator}}\\
\vspace{+5mm}
\normalsize{\bf{Celal Umut Yaran$^{\dag}$*\footnote{\textbf{Corresponding author}, cyaran@itu.edu.tr}, Mine Çağlar*}}
\end{center}
    \normalsize{$\dag$ Department of Mathematics Engineering, İstanbul Technical University }\\
     \normalsize{* Mathematics Department, Koç University}

\begin{abstract}
\noindent
   Our main result is the martingale representations for Markov additive processes where the modulator is a Lévy process. These processes have three parts: the modulator, the jumps of the ordinate triggered by the modulator, and the semimartingale part of the ordinate with parameters depending on the modulator. We orthogonalize Teugels martingales constructed from these parts to give a chaotic representation of square-integrable random variables as a sum of stochastic integrals with respect to the orthogonal sequence obtained. Consequently, a predictable representation of square-integrable martingales is derived in terms of the ordinate and the Teugels martingales.
\end{abstract}

\section{Introduction}
The classical Itô theorem establishes that every square-integrable martingale adapted to the natural filtration of a Brownian motion $B=(B_t)_{t\geq 0}$ can be represented as a stochastic integral with respect to the same Brownian motion \cite[Thm.36.1]{MR1780932}. Formally, it states that if $M$ is a square-integrable martingale adapted to the natural filtration $\mathcal{F}=(\mathcal{F}_t)_{t\geq 0}$ of $B$, then there exists a predictable process $h_M$ adapted to $\mathcal{F}$ such that\begin{equation*}
    M_t=M_0 + \int_0^t h_M(t)dB_t.
\end{equation*} A key application of martingale representations in mathematical finance is determining whether a financial market is complete by calculating the hedging strategy \cite[$\S$1.6]{MR1640352}. In this case, it is shown that each claim can be replicated by a self-financing portfolio under certain assumptions \cite[Cor.1.6.8]{MR1640352}.

The classical martingale representation theorem has further generalizations \cite{MR217856,MR2512800,MR418221,MR420846,MR433584}. In \cite{MR1787127}, chaotic and predictable representation theorems for Lévy processes has been derived. It is shown that all square-integrable martingales adapted to the natural filtration generated by a Lévy process $\xi$ can be written as a sum of stochastic integrals with respect to certain strongly orthogonal martingales constructed from $\xi$. An application to a Lévy market model is given in \cite{MR2210930}. 

Another generalization is the martingale representation property for a Markov additive process(MAP). A MAP is a natural extension of Lévy processes and they were introduced in \cite{MR329047}, as a real-valued stochastic process $(\xi_t)_{t\geq 0}$ such that the increments of $\xi$ are modulated by a Markov process $(\Theta_t)_{t\geq 0}$. Formally, the process $(\xi, \Theta)$ is called a MAP if given $\{ (\xi_s,\Theta_s),s\leq t\}$ for any $t\geq 0$, the process $(\xi_{s+t}-\xi_t, \Theta_{s+t})_{s\geq 0}$ has the same law as $(\xi_s,\Theta_s)_{s\geq 0}$ under $\mathbb{P}_{0,v}$ with $v=\Theta_t$, where $\mathbb{P}_{x,\theta}(\xi_0=x,\Theta_0=\theta)=1$. The following studies on general Markov additive processes can be found first in \cite{MR370788,MR445621,MR650610,MR714146}, and in \cite{MR885131,MR901545}. Although MAPs have been introduced in such generality, most of the previous work has focused on the case the \textit{modulator} $\Theta$ is a Markov chain with a discrete state space (e.g. \cite{MR3866614,MR2766220,MR478365}). In this case, the chaotic and predictable representation properties have been derived in \cite{MR3854533}. In this paper, we generalize this result for a MAP with a Lévy modulator. The recent studies on MAPs when the state space of the modulator is a Polish space can be found in \cite{MR3622890,MR3808900,MR4692990,MR4155177,MR4801600,MR3912200,MR4958275}.

Let $(\xi,\Theta)$ be a MAP on a filtered probability space $(\Omega,\mathcal{F}_\infty,(\mathcal{F}_t)_{t\geq 0}, \mathbb{P})$  such that the modulator $\Theta$ is a real Lévy process. Then, we can decompose the ordinate $\xi$ as\begin{equation*}
    \xi = \xi^L + \xi^f
\end{equation*}
where $\xi^L$ is a real-valued semimartingale such that its jumps are independent of $\Theta$, and $\xi^f$ is a pure jump process such that its jump times are fixed by $\Theta$, and its jump size depends on those of $\Theta$. Indeed, it is a particular form of the decomposition of the ordinate in the general case given in \cite{MR329047}. We introduce compensated Teugels jump processes  $\bar{\Theta}^{(l)}$, $\bar{\xi}^{(k)}$, and $\bar{\xi}^{[m]}$ for $k,l,m\geq 1$, corresponding to $\Theta$, $\xi^L$, and $\xi^f$, respectively. These additional processes are needed for a martingale representation. Then, we construct an orthogonal sequence of martingales $\{H^{(k)},k\geq 1\}$ obtained by orthogonalization of $\bar{\xi}^{(k)}$, and another sequence $\{ G^{(k,l)},k\geq 1, l\geq 0\}$ found by orthogonalization of $\bar{\Theta}^{(l)}$ and $\bar{\xi}^{[m]}$.

In Proposition \ref{Representation of Powers}, we show that monomials of the form $\bar{\Theta}^g (\bar{\xi}^L)^p(\bar{\xi}^f)^b$ for $g,p,b\geq 0$ can be represented as a sum of stochastic integrals with respect to the integrators $\bar{\xi}^{(k)}$, $\bar{\Theta}^l$, and $\bar{\xi}^{(m)}$ for $k,l,m\geq 1$, where $\bar{\Theta}$, $ \bar{\xi}^L$, and $\bar{\xi}^f$ are compensated versions of $\Theta$, $\xi^L$, and $\xi^f$, respectively. Then, in Proposition \ref{Base}, we show that $L^2(\Omega,\mathcal{F}_t)$ is generated by monomials of the form  $\{ \bar{\Theta}_{t_1}^{g_1}\cdot \ldots \cdot \bar{\Theta}_{t_m}^{g_m}\cdot (\bar{\xi}_{t_1}^L)^{p_1}\cdot \ldots \cdot (\bar{\xi}_{t_m}^L)^{p_m}\cdot (\bar{\xi}_{t_1}^f)^{b_1}\cdot \ldots \cdot (\bar{\xi}_{t_m}^f)^{b_m} \, : 0\leq t_1<\ldots <t_m\leq t\}$ for $g_1,\ldots,g_m,p_1,\ldots,p_m,b_1,\ldots,b_m\geq 0$ 
under certain conditions. The main results of the paper are the representation properties for a MAP $(\xi,\Theta)$. The chaotic representation property given in Theorem \ref{Chaotic Representation Theorem} states that any random variable $F$ in $L^2(\Omega,\mathcal{F}_t)$ can be represented as\begin{align*}
        F = & \mathbb{E}[F] + \sum_{s=1}^\infty \sum_{\tau=1}^\infty \sum_{\zeta=1}^\infty \sum_{\nu_1,\ldots,\nu_\zeta\geq 1} \sum_{\mu_1,\ldots,\mu_\tau\geq 0 } \sum_{\iota_1,\ldots,\iota_s\geq 1}
         \int_0^t\int_0^{t_{1^-}}\ldots \int_0^{t_{\zeta+s+\tau-1}^-}\\
         &f_{\nu_1,\ldots,\nu_\zeta,\mu_1,\ldots,\mu_\tau,\iota_1,\ldots,\iota_s}(t,t_1,t_2,\ldots,t_{\zeta+s+\tau})
         dH_{t_\zeta+s+\tau}^{(\nu_\zeta)}\ldots dH_{t_{1+s+\tau}}^{(\nu_1)}dG_{t_{s+\tau}}^{(\iota_s,\mu_\tau)}\ldots dG_{t_1}^{(\iota_1,\mu_0)}
    \end{align*}
   where $f_{(\nu_1,\ldots , \nu_\zeta,\mu_1,\ldots , \mu_\tau, \iota_1, \ldots, \iota_s)}(t,t_1,\ldots, t_{\zeta+s+\tau})$ are some random fields such that the integral is well-defined on $L^2(\Omega,\mathcal{F}_t)$. In Theorem \ref{Main Predictable Representation Theorem}, we give the representation of any square integrable $(\mathcal{F}_t)$-martingale $M$ as a sum of stochastic integrals of predictable processes with respect to $\xi$ and the Teugels martingales $\bar{\xi}^{(l)}$, $\bar{\Theta}^{(l)}$ and $\bar{\xi}^{[m]}$ as \begin{equation*}
        M_t=M_0 +\int_0^t h_\xi(s)d\bar{\xi}_s + \sum_{k=2}^\infty \int_0^t h_L^{(k)}(s)(s)d\bar{\xi}_s^{(k)}+ \sum_{l=1}^\infty \int_0^t h_f^{(l)}(s)d\bar{\xi}_s^{[l]} 
        +\sum_{m=1}^\infty \int_0^t h_\Theta^{(m)}(s)d\bar{\Theta}_s^{(m)}
    \end{equation*}
where $h_\xi,h_L^{(k)},h_f^{(l)}$, and $h_\Theta^{(m)}$ are predictable processes.
   This property is known as the predictable representation property.

The paper is organized as follows. In Section \ref{Preliminaries}, we define our MAP through a Lévy modulator and an ordinate with a decomposition into a semimartingale and a pure-jump process. The Teugels power jump processes and the corresponding orthogonal sequence of martingales are constructed in Section \ref{Section: Orthogonal Martingales}. Finally, the main results for the paper, namely chaotic and predictable representation properties for a MAP, are presented in Section \ref{Section: Representations}.

\section{A Markov Additive Process with Lévy Modulator}\label{Preliminaries}
Let $(\xi_t,\Theta_t)_{t\geq 0}$ be a possibly killed Markov process on a filtered probability space \linebreak
$(\Omega,\mathcal{F}_\infty,(\mathcal{F}_t)_{t\geq 0},\mathbb{P}_{x,\theta})$ with $\mathbb{P}_{x,\theta}(\xi_0=x,\Theta_0=\theta)=1$ and state space $(\mathbb{R}\times S,\mathcal{B}(\mathbb{R}\times S))$ with an extra isolated state $\partial$, where $S$ is a Polish space, $\mathcal{B}(\mathbb{R}\times S)$ denotes the Borel $\sigma$-algebra on $\mathbb{R}\times S$, $(\mathcal{F}_t)_{t\geq 0}$ is the minimal augmented admissible filtration and $\mathcal{F}_\infty = \vee_{t=0}^\infty \mathcal{F}_t$.
The process $(\xi_t,\Theta_t)_{t\geq 0}$ is called a Markov additive process (MAP) on $\mathbb{R}\times S$, if \begin{equation*}
     \mathbb{P}_{0,\theta}( (\xi_{t+s}-\xi_t)\in \Gamma \, , \, \Theta_{t+s}\in A \, | \, \mathcal{F}_t ) = \mathbb{P}_{0,v}(\xi_s \in \Gamma \, , \, \Theta_s\in A)
    \end{equation*}
for all $\theta\in S$, $\Gamma\in \mathcal{B}( \mathbb{R} )$, $A \in \mathcal{B} (S)$, and $t\geq 0$,  where $v=\Theta_t$, and denoted by $((\xi,\Theta),\mathbb{P})=((\xi_t,\Theta_t)_{t\geq 0},\mathbb{P}_{0,\theta})$ \cite{MR329047}. We call $\xi$ as the ordinate and $\Theta$ as the modulator of the MAP.

Consider a Lévy process $\Theta$ generating the $\sigma$-algebra, $\mathcal{K}_t=\sigma(\Theta_s\, , \, s\leq t)$, $t\geq 0$, and the following representation\begin{equation*}
        \Theta_t= \mu_1 t + \sigma_1 B_t   + \int_0^t \int_\mathbb{R}  x\bar{M_1}(ds,dx)
    \end{equation*}
    where $\mu_1,\sigma_1\in \mathbb{R}\backslash\{0\}$ are constants, $B$ is a real-valued Brownian motion, $M_1$ is a Poisson random measure on $\mathbb{R}\times \mathbb{R}$ with mean measure $\text{Leb}\times \nu_1$ such that \begin{equation*}
        \int_\mathbb{R}\nu_1(dx)\min{(|x|^2,1)}<\infty
    \end{equation*}
    and \begin{equation*}
         \bar{M}_1(dx,ds):= \begin{cases}
         M_1(ds,dx)-\nu_1(dx)ds & \text{if  }  |x|<1  \\
        M_1(ds,dx) &  \text{if   }  |x|\geq 1.
    \end{cases}
    \end{equation*}
Let $\xi$ be a real-valued quasi-left continuous process with the following decomposition\begin{equation}\label{decomposition}
    \xi_t = \xi_t^L + \xi_t^f
\end{equation}
where $\xi_t^L$ is a semimartingale with lifetime $T\geq 0$ and has the following representation for $t\leq T$:\begin{equation*}
    \xi_t^L = \xi_0 + \int_0^t \mu_2(\Theta_s)ds + \int_0^t \sigma_2(\Theta_s)dW_s + \int_0^t\int_\mathbb{R} x \bar{M}_2(ds,dx)
\end{equation*}
    where $\mu_2,\sigma_2$ are real-valued functions such that almost surely\begin{equation*}
        \int_0^T \sigma_2^2(\Theta_s)ds <\infty \quad , \quad \int_0^T|\mu_2(\Theta_s)|ds<\infty ,
    \end{equation*} $W$ is a real-valued Brownian motion independent of $\Theta$, and $M_2$ is a Poisson random measure on $\mathbb{R}\times \mathbb{R}$ independent of $\Theta$ and $W$ with mean measure $\text{Leb}\times \nu_2$ such that\begin{equation*}
        \int_\mathbb{R}\nu_2(dx)\min{(|x|^2,1)}<\infty
    \end{equation*}
    and \begin{equation*}
         \bar{M}_2(dx,ds):= \begin{cases}
         M_2(ds,dx)-\nu_2(dx)ds & \text{if  }  |x|<1  \\
        M_2(ds,dx) &  \text{if   }  |x|\geq 1.
    \end{cases}
    \end{equation*}
    Note that $\xi_t^L$ is a jump-diffusion process given $\{ \Theta_s, s\leq t\}$ whose the diffusion and drift parameters depend on the state of the $\Theta$.
    On the other hand, $\xi_t^f$ is a pure-jump process of the form\begin{equation*}
        \xi_t^f = \sum_{j\geq 1} U_j(\Delta\Theta_{\tau_j})\mathbf{1}_{\tau_j \leq t}
    \end{equation*}
    where $\tau_j$'s are jump times of $\Theta_t$, $\Delta(\Theta_t):= \Theta_t - \Theta_{t^-}$, and $U_j$'s are independent identical pure jump processes. Note that $\xi_t^L$ and $\xi_t^f$ are conditionally independent given $\mathcal{K}_t$. We can express $\xi_t^f$ as\begin{equation*}
        \xi_t^f = \int_0^t \int_\mathbb{R} \int_\mathbb{R} y \Pi_U(ds,dx,dy)
    \end{equation*}
    where \begin{equation*}
        \Pi_U([0,t],dx,dy)= \sum_{j\geq 1}\mathbf{1}_{\{U_j(x)\in dy\}} \mathbf{1}_{\{ \Delta \Theta_{\tau_j}\in dx,\tau_j \leq t \}}
    \end{equation*}
Then, $(\xi,\Theta)$ is a $\mathbb{R}\times\mathbb{R}$-valued MAP. 
\begin{Rem}
    Although $(\xi,\Theta)$ is a special MAP, the decomposition \eqref{decomposition} of $\xi$ is valid in the general case given in \cite{MR329047}. The ordinate $\xi$ is decomposed as \begin{equation*}
    \xi = A + \xi^c + \xi^d + \xi^f
\end{equation*}
where the first three terms are similar to the decomposition of a Lévy process; that is, $A$ is an additive functional of $\xi$, $\xi^c$ is a continuous process, and $\xi^d$ is a pure-jump process whose jump times are independent of the modulator $\Theta$. On the other hand, $\xi^f$ is a pure-jump process such that jump times are fixed by $\Theta$. For example, when $\Theta$ is a finite state space Markov chain, then a transition of $\Theta$ triggers a jump of $\xi$ which is represented in the part $\xi^f$ \cite{MR3854533}. In general, according to \cite[Thm.4.8]{MR329047}, the conditional distribution $U_j$ of the jumps given $\mathcal{K}_{\tau_j}$ only depends on $\Theta_{\tau_j^-}$ and $\Theta_{\tau_j}$ under certain assumptions. In our framework, we use the notation $\xi^L$ for the Lévy-like part $\xi^L=A+\xi^c+\xi^d$ and $\xi^f$ for the part that consists of jumps triggered from the modulator. Furthermore, as a particular case, we assume that the conditional distribution of the triggered jumps depend on the jump size of the modulator, that is, $\Theta_{\tau_j}-\Theta_{\tau_j^-}$.
\end{Rem}
We assume that for some $\epsilon,\lambda >0$ \begin{itemize}
    \item $\int_{|x|\geq 1}e^{\lambda x}\nu_i(dx)<\infty$ a.s. for $i=1,2$.
    \item $\int_{|x|\geq 1}e^{\lambda x} \mathbb{P}(U(y)\in dx) <\infty$ a.s. for all $y\in \mathbb{R}$
\end{itemize}
which imply that  $\xi$ has moments of all orders and polynomials are dense in $L^2(\mathbb{R},\varphi_t)$ where $\varphi_t$ is the probability distribution of $\xi_t$ \cite[$\S$2]{MR3854533}. 

\section{Power Jump Processes and Orthogonal Martingales}\label{Section: Orthogonal Martingales}

The chaotic and predictable representation properties for a Lévy process have been derived in \cite{MR1787127}. These representations are given in terms of an orthogonal sequence of martingales obtained as transformations of power jump processes which are associated with the Lévy process. This idea is extended to the representation theorems for an Itô-Markov additive process in \cite{MR3854533} ,where the modulator is a finite-state continuous time Markov chain. Our main aim is to extend these results for a MAP defined in Section \ref{Preliminaries} where the modulator is a Lévy process. In this section, first we define the power jump processes and then derive the orthogonal martingales.

We work with augmented filtration $\mathcal{F}_t=\sigma(\Theta_s,\xi_s \, : s\leq t)\vee \mathcal{N}$ where $\mathcal{N}$ denotes the $\mathbb{P}$-null sets of $\mathcal{F}$. This filtration is equaled to \begin{equation}\label{sigma algebra}
    \mathcal{F}_t = \sigma (\Theta_s, \xi_s^L,\xi_s^f \, : \, s\leq t)\vee \mathcal{N}
\end{equation}
by the same arguments in the proof of Theorem 3.3 of \cite{MR2060091} as in \cite{MR3854533}. We now construct the power jump processes associated with $\Theta$, $\xi^L$ and $\xi^f$ as analogous to \cite{MR3854533,MR1787127}. First, define the following processes associated with $\xi^L$ \begin{equation*}
    \xi_t^{(k)}= \begin{cases}
        \sum_{0<s\leq t}(\Delta \xi_s^L)^k & \text{for  }  k \geq 2  \\
        \xi_t^L &  \text{for   }  k=1.
    \end{cases}
\end{equation*}
Each process $\xi^{(k)}$ has the same jump times as the original process $\xi^L$, but the sizes of the jumps are $k^{th}$ power of the jump sizes of $\xi^L$. Furthermore, they are adapted to $(\mathcal{F}_t)$ \cite{MR1943877,MR3854533}. By \cite{MR2273672}, one can show that
    $\mathbb{E}[\xi_t^{(k)}\,|\,\mathcal{K}_t]= \int_0^t \int_\mathbb{R} x^k\nu_2(dx)ds$ for $k\geq 2$. Then, the following are $(\mathcal{F}_t)$-martingales:\begin{equation*}
    \bar{\xi_t}^{(k)} = \begin{cases}
        \xi_t^{(k)}-\int_0^t \int_\mathbb{R} x^k\nu_2(dx)ds & \text{for  }  k \geq 2  \\
        \xi_t^L - \mathbb{E}[ \xi_t^L \, | \, \mathcal{K}_t] &  \text{for   }  k=1,
    \end{cases}
\end{equation*}
called Teugels martingales \cite{MR3854533,MR1761401}. 

Similarly, we can define the power jump processes associated with the modulator $\Theta$ as follows
\begin{equation*}
    \Theta_t^{(k)}= \begin{cases}
        \sum_{0<s\leq t}(\Delta \Theta_s)^k & \text{for  }  k \geq 2  \\
        \Theta_t &  \text{for   }  k=1
    \end{cases}
\end{equation*}
Then, one can show that  
    $\mathbb{E}[\Theta_t^{(k)}]= \int_0^t \int_\mathbb{R} x^k\nu_1(dx)ds$ for $k\geq 2$ to obtain the following $(\mathcal{F}_t)$-martingales:\begin{equation*}
    \bar{\Theta}_t^{(k)} = \begin{cases}
        \Theta_t^{(k)}-\int_0^t \int_\mathbb{R} x^k\nu_1(dx)ds & \text{for  }  k \geq 2  \\
        \Theta_t - \mathbb{E}[ \Theta_t ] &  \text{for   }  k=1
    \end{cases}
\end{equation*}

Finally, we define the corresponding power jump processes associated with $\xi_t^f$\begin{equation*}
    \xi_t^{[k]}= \int_0^t\int_\mathbb{R}\int_\mathbb{R} y^k\Pi_U(ds,dx,dy).
\end{equation*}
Then, the processes\begin{equation*}
    \bar{\xi}_t^{[k]}= \int_0^t\int_\mathbb{R}\int_\mathbb{R}y^k\Tilde{\Pi}_U(ds,dx,dy)
\end{equation*}
are $(\mathcal{F}_t)$-martingales where $\Tilde{\Pi}_U(ds,dx,dy):= \Pi_U(ds,dx,dy)-\nu_1(dx)\mathbb{P}(U(x)\in dy)ds$. For the underlying motivations for these processes, we refer to \cite{MR3854533,MR2210930}.

Let $\mathcal{M}^2$ denote the set of square-integrable $(\mathcal{F}_t)$-martingales. Two martingales $M_1,M_2\in \mathcal{M}^2$ are called to be strongly orthogonal if their product $M_1 \cdot M_2$ is a uniformly integrable martingale. $M_1,M_2 \in \mathcal{M}^2$ are strongly orthogonal if and only if their cross variation $[M_1,M_2]$ is a uniformly integrable martingale \cite{MR2273672}. We need to construct a set of pairwise strongly orthogonal martingales in the form of linear combinations of $\bar{\xi}^{(k)}$, $\Theta^{(l)}$, and $\bar{\xi}^{[m]}$. We start with obtaining $\{ H^{(k)}, k\geq 1 \}$ of pairwise strongly orthogonal martingales of the form \begin{equation}\label{H}
    H^{(k)}= a_{k,k}\bar{\xi}^{(k)} + a_{k,k-1}\bar{\xi}^{(k-1)} + \ldots + a_{k,1}\bar{\xi}^{(1)}
\end{equation}
where $a_{k,l}\in\mathbb{R}$ for all $k,l\geq 1$. We use the orthogonalization of polynomials to compute the coefficients $a_{k,l}$ as in \cite{MR3854533,MR1787127}.  Let $S_1$ be the space of all real-valued polynomials and define an inner product $\langle\, , \rangle$ for some $\alpha>0$ as\begin{equation*}
    \langle P(x),Q(x)\rangle_1 = \int_\mathbb{R}P(x)Q(x)x^2\nu_2(dx) + \alpha P(0)Q(0)
\end{equation*}
Let $S_2$ denote the space of all linear transformations of $\bar{\xi}^{(k)}$ for $k\geq 1$ with the inner product $\langle\, , \rangle_2$ defined as\begin{equation*}
    \langle\bar{\xi}^{(k)},\bar{\xi}^{(l)}\rangle_2 = \mathbb{E}[\,[\bar{\xi}^{(k)}_1,\bar{\xi}^{(l)}_1]\,]
\end{equation*}
Note that \begin{equation*}
[\bar{\xi}^{(k)}_1,\bar{\xi}^{(l)}_1]= \int_\mathbb{R}x^{k+l}\nu_2(dx) + \int_0^1\sigma_2^2(\Theta_s)ds\mathbf{1}_{\{k=l=1\}}
\end{equation*}
Then, by choosing $\alpha =\int_0^1\sigma_2^2(\Theta_s)ds $, the correspondence $x^{k-1} \leftrightarrow \bar{\xi}^{(k)} $ for $k\geq 1$ gives an isometry between the spaces $S_1$ and $S_2$. Hence an orthogonalization of $\{1,x,x^2,\ldots \}$ on $S_1$ gives an orthogonalization of $\{ \bar{\xi}^{(1)},\bar{\xi}^{(2)},\ldots \}$ in $L^2(\Omega,\mathcal{F})$.\\

Similarly, we want to obtain $\{ G^{(k,l)},k\geq 1,l\geq 0\}$ of pairwise strongly orthogonal martingales such that\begin{equation}\label{G}
    G^{(k,l)}= b_{k,l}\bar{\Theta}^{(k)} + b_{k-1,l}\bar{\Theta}^{(k-1)}+\ldots +b_{1,l}\bar{\Theta}^{(1)}+ b_{k,l-1}\bar{\xi}^{[l]}+b_{k,l-1}\bar{\xi}^{[l-2]}+\ldots + b_{k,0}\bar{\xi}^{[1]}
\end{equation}
Define the following inner product space $(S_3,\langle\, , \rangle_3)$\begin{equation*}
    S_3 = \{ P(x,y)=f(x)+g(y) \, | \, f\text{ and }g \text{ are real-valued polynomials }, \, g(0)=0 \} 
    \end{equation*}
    with the following inner product,
    \begin{align*}
    \langle P_1(x,y),P_2(x,y)\rangle_3 &=\int_\mathbb{R}\int_\mathbb{R} g_1(y)g_2(y)\nu_1(dx)\mathbb{P}(U(x)\in dy) +\\
    &\int_\mathbb{R}\int_\mathbb{R}\left(g_1(y)f_2(x)+g_2(y)f_1(x)\right)x\nu_1(dx)\mathbb{P}(U(x)\in dy) +\\
    &\int_\mathbb{R}f_1(x)f_2(x)x^2\nu_1(dx) + \sigma_1^2f_1(0)f_2(0) 
\end{align*}
where $P_1(x,y)=f_1(x)+g_1(y)$, and $P_2(x,y)=f_2(x)+g_2(y)$. It can be deduced that $(S_3,\langle\, , \rangle_3)$ is an inner product space. Actually, linearity and symmetry is straightforward. We show the positive-definiteness. Let $P(x,y)=f(x)+g(y)\neq 0$. Then,\begin{align*}
    <P(x,y),P(x,y)>_3 = &\int_\mathbb{R}\int_\mathbb{R} g^2(y)\nu_1(dx)\mathbb{P}(U(x)\in dy) \\
    &+\int_\mathbb{R}\int_\mathbb{R}2f(x)g(y)x\nu_1(dx)\mathbb{P}(U(x)\in dy) \\
    &+\int_\mathbb{R}f^2(x)x^2\nu_1(dx) + \sigma_1^2f^2(0)  \\
    = &\int_\mathbb{R}\int_\mathbb{R}\left(g^2(y)+ 2f(x)g(y)x + f^2(x)x^2\right)\nu_1(dx)\mathbb{P}(U(x)\in dy)\\
    &+ \sigma_1^2f^2(0) \\
    =&\int_\mathbb{R}\int_\mathbb{R}(g(y) +xf(x))^2\nu_1(dx)\mathbb{P}(U(x)\in dy)+\sigma_1^2f^2(0) 
\end{align*}
If $g(y)+xf(x)\neq 0$, this integral is positive. Assume that $g(y)+xf(x)=0$ identically. Then, since it holds for each pair of the form $(0,y)$, it implies that $g(y)=0$ identically. We obtain that $xf(x)=0$ for all $x$, which implies that $f(x)=0$ for all $x\neq 0$. In this case, $P(x,y)= f(x) + g(y)= f(x)=0$ for all pairs $(x,y)$ with $x\neq0$. Since $P(x,y)$ is not identically equals $0$, it implies that $f(0)\neq 0$, so $\sigma^2_1f^2(0)> 0$ and we are done.\\

Let $S_4$ denote the space of all linear transformations of the form \eqref{G} with an inner product space defined as\begin{eqnarray*}
    \langle a,b\rangle_4 = \mathbb{E}[\,[a_1,b_1]\,]
\end{eqnarray*}
In particular, we have\begin{align*}
    \langle\bar{\xi}^{[k]},\bar{\xi}^{[l]}\rangle_4 &= \int_\mathbb{R}\int_\mathbb{R} y^{k+l}\nu_1(dx)\mathbb{P}(U(x)\in dy)\\
    \langle\bar{\xi}^{[k]},\bar{\Theta}^{(l)}\rangle_4 &= \int_\mathbb{R}\int_\mathbb{R} y^{k}x^{l} \nu_1(dx)\mathbb{P}(U(x)\in dy)\\
    \langle\bar{\Theta}^{(k)},\bar{\Theta}^{(l)}\rangle_4 &= \int_\mathbb{R}x^{k+l}\nu_1(dx) + \sigma^2 \mathbf{1}_{\{k=l=1\}}
\end{align*}
Then, the correspondence \begin{align*}
    x^{i-1} &\leftrightarrow \bar{\Theta}^{(i)} \\
    y^i & \leftrightarrow \bar{\xi}^{[i]}
\end{align*}
gives an isometry between $S_3$ and $S_4$. An orthogonalization of  $\{ 1,x,x^2, \ldots, x^{k-1}$ $,y,y^2,\ldots, y^l\}$ in $S_3$ gives an orthogonalization of $\{ \bar{\Theta}^{(1)}, \ldots \bar{\Theta}^{(k)}, \bar{\xi}^{[1]},\ldots, \bar{\xi}^{[l]}\}$. Since $\xi^L$ does not jump at the same time with $\xi^f$ and $\Theta$, and the continuous parts of $\xi^L$ and $\Theta$ are independent, the martingales $H^{(k)}$ and $G^{(l,m)}$ are pairwise strongly orthogonal. Indeed, $[H^{(k)},G^{(l,m)}]$ equals $0$ and trivially is a uniformly integrable martingale. Hence, they are pairwise strongly orthogonal.

\section{Chaotic and Predictable Representation Properties}\label{Section: Representations}

We are ready to give chaotic and predictable representation theorems for the MAP $(\xi,\Theta)$. Let $\bar{\xi}_t^L:= \xi_t^L - \mathbb{E}[\xi_t^L | \mathcal{F}_t] =\bar{\xi}_t^{(1)}$, $\bar{\Theta}_t:=\Theta_t - \mathbb{E}[\Theta_t]=\bar{\Theta}_t^{(1)}$, and $\bar{\xi}_t^f:= \int_0^t \int_\mathbb{R}\int_\mathbb{R}y\tilde{\Pi}_U(ds,dx,dy)=\bar{\xi}_t^{[1]}$ denote the compensated versions of $\xi^L,\Theta,\xi^f$, respectively. As in \cite{MR3854533, MR1787127}, we first give a representation of monomials of the form $\bar{\Theta}_t^g (\bar{\xi}_t^L)^p (\bar{\xi}_t^f)^b$ in terms of the martingales $\bar{\xi}^{(k)},\bar{\Theta}^{(l)},\bar{\xi}^{[m]}$ in the next proposition. 

\begin{Prop}\label{Representation of Powers}
    For $g,d,p\geq 0$, we have the representation\begin{align*}
        \bar{\Theta}_t^g (\bar{\xi}_t^L)^p (\bar{\xi}_t^f)^b =& f^{(g+p+b)}(t)+\sum_{\zeta=1}^p \sum_{s=1}^g\sum_{\tau=1}^b\sum_{(\nu_1,\ldots,\nu_\zeta)\in \{1,\ldots,p\}^\zeta} \sum_{(\mu_1,\ldots,\mu_\tau)\in \{1,\ldots,b\}^\tau}\sum_{(\iota_1,\ldots,\iota_s)\in \{1,\ldots,g\}^s}\\
        & \int_0^t \int_0^{t_1^-} \ldots \int_0^{t_{(\zeta + s + \tau -1)}^-} f_{(\nu_1,\ldots , \nu_\zeta,\mu_1,\ldots , \mu_\tau, \iota_1, \ldots, \iota_s)}^{(g+p+b)}(t,t_1,\ldots, t_{\zeta+s+\tau}) \\ & \quad \quad \quad \quad \quad \quad \quad \quad d\bar{\xi}_{t_{\zeta+s+\tau}}^{(\nu_\zeta)} \ldots d\bar{\xi}_{t_{1+s+\tau}}^{(\nu_1)}d\bar{\Theta}_{t_{s+\tau}}^{(\iota_s)}\ldots d\bar{\Theta}_{t_{\tau+1}}^{(\iota_1)}d\bar{\xi}_{t_{\tau}}^{[\mu_\tau]}\ldots d\bar{\xi}_{t_1}^{[\mu_1]}
    \end{align*} where $f^{(g+p+b)}(t)$ and $f_{(\nu_1,\ldots , \nu_\zeta,\mu_1,\ldots , \mu_\tau, \iota_1, \ldots, \iota_s)}^{(g+p+b)}(t,t_1,\ldots, t_{\zeta+s+\tau})$ are sums of predictable processes defined on $L^2(\Omega,\mathcal{F}_t)$.
\end{Prop}
\textbf{Proof:} We follow the argument of \cite[Thm.7]{MR3854533}. We prove this by induction. Indeed, we 
represent the product $\bar{\Theta}^g (\bar{\xi}^L)^p (\bar{\xi}^f)^b$ as a sum of stochastic integrals of 
lower powers of $\bar{\Theta}$, $\bar{\xi}^L$, $\bar{\xi}^f$ with respect to processes $\bar{\Theta}^{(k)}$, 
$\bar{\xi}^{(l)}$, and $\bar{\xi}^{[m]}$ for $k\leq g$, $l\leq p$ and $m\leq b$. Note that 
$[\bar{\Theta}_t,\bar{\Theta}_t]^c= \sigma_1^2t$, 
$[\bar{\xi}^L_t,\bar{\xi}^L_t]^c=\int_0^t\sigma_2^2(\Theta_{s^-})ds$, and $[\bar{\xi}^f_t,\bar{\xi}^f_t]^c=0$. 
Furthermore, by \cite[Thm.4.19]{MR2933773}, and \cite[pg.70]{MR2273672}, we have 
$[\bar{\Theta}_t,\bar{\xi}^L_t]^c=[\bar{\Theta}_t,\bar{\xi}^f_t]^c = [\bar{\xi}^L_t,\bar{\xi}^f_t]^c=0$. By Itô's Lemma, we obtain \begin{align}
    \bar{\Theta}_t^g (\bar{\xi}_t^L)^p (\bar{\xi}_t^f)^b &= 
    \bar{\Theta}_0^g(\bar{\xi}_0^L)^p(\bar{\xi}_0^f)^b + \int_{0^+}^t g\bar{\Theta}_{s^-}^{g-1}(\bar{\xi}_{s^-
    }^L)^p(\bar{\xi}_{s^-}^f)^b d\bar{\Theta}_s + \int_{0^+}^tp\bar{\Theta}_{s^-}^g (\bar{\xi}_{s^-}^L)^{p-1}
    (\bar{\xi}_{s^-}^f)^bd\bar{\xi}_s^L + \nonumber  \\
    & \int_{0^+}^t b\bar{\Theta}_{s^-}^g(\bar{\xi}_{s^-}^L)^p(\bar{\xi}_{s^-}^f)^{b-1}d\bar{\xi}_{s}^f + \frac{1}{2}p(p-1) \mathbf{I_1} + \frac{1}{2}g(g-1)\mathbf{I_2} + \mathbf{I_3} \label{1.}
\end{align}
where \begin{equation*}
    \mathbf{I_1}= \int_{0^+}^t \bar{\Theta}_{s^-}^g(\bar{\xi}_{s^-}^L)^{p-2}(\bar{\xi}_{s^-}^f)^b\sigma_2^2(\bar{\Theta}_{s^-})ds
\end{equation*}
\begin{equation*}
    \mathbf{I_2}= \int_{0^+}^t \bar{\Theta}_{s^-}^{g-2}(\bar{\xi}_{s^-}^L)^{p}(\bar{\xi}_{s^-}^f)^b\sigma_1^2ds
\end{equation*}
\begin{align*}
\mathbf{I_3}= \sum_{0<s\leq t} \left[\ 
   \bar{\Theta}_s^g(\bar{\xi}_s^L)^p(\bar{\xi}_s^f)^b - \bar{\Theta}_{s^-}^g(\bar{\xi}_{s^-}^L)^p(\bar{\xi}_{s^-}^f)^b-g\bar{\Theta}_{s^-}^{g-1}(\bar{\xi}_{s^-}^L)^p(\bar{\xi}_{s^-}^f)^b\Delta\bar{\Theta}_s \right. \\
   \left. -p\bar{\Theta}_{s^-}^{g}(\bar{\xi}_{s^-}^L)^{p-1}(\bar{\xi}_{s^-}^f)^b\Delta\bar{\xi}_s^L-
   b\bar{\Theta}_{s^-}^{g}(\bar{\xi}_{s^-}^L)^p(\bar{\xi}_{s^-}^f)^{b-1}\Delta\bar{\xi}_s^f \right].   
\end{align*}
For $\mathbf{I_1}$, we apply integration by parts for the process $\bar{\Theta}_t^g(\bar{\xi}_t^L)^{p-2}(\bar{\xi}_t^f)^b\int_0^t\sigma^2_2(\bar{\Theta}_s)ds$. Then, we obtain\begin{align*}
    \bar{\Theta}_t^g(\bar{\xi}_t^L)^{p-2}(\bar{\xi}_t^f)^b\int_0^t\sigma^2_2(\bar{\Theta}_s)ds &= \int_0^t \bar{\Theta}_{s^-}^g(\bar{\xi}_{s^-}^L)^{p-2}(\bar{\xi}_{s^-}^f)^b\sigma^2_2(\bar{\Theta}_{s^-})ds  \\
    &+ \int_0^t\left(\int_0^s \sigma_2^2(\bar{\Theta}_{u^-})du\right) d( \bar{\Theta}_s^g(\bar{\xi}_s^L)^{p-2}(\bar{\xi}_s^f)^b) \\
    &+  [ \bar{\Theta}_t^g(\bar{\xi}_t^L)^{p-2}(\bar{\xi}_t^f)^b,\int_0^t\sigma^2_2(\bar{\Theta}_s)ds] 
\end{align*}
The given quadratic covariation in the last term is zero. Consider the differential $d( \bar{\Theta}_s^g(\bar{\xi}_s^L)^{p-2}(\bar{\xi}_s^f)^b)$ in the second term. Using integration by parts, we obtain\begin{eqnarray*}
    d( \bar{\Theta}_s^g(\bar{\xi}_s^L)^{p-2}(\bar{\xi}_s^f)^b)= (\bar{\xi}_{s^-}^L)^{p-2}d(\bar{\Theta}_s^g(\bar{\xi}_s^f)^b)+ \bar{\Theta}_{s^-}^g(\bar{\xi}_{s^-}^f)^bd(\bar{\xi}_s^L)^{p-2}+d[(\bar{\xi}_s^L)^{p-2}, \bar{\Theta}_s^g(\bar{\xi}_s^f)^b)]
\end{eqnarray*}
Since the continuous parts of the processes $(\bar{\xi}_s^L)^{p-2}$ and $\bar{\Theta}_s^g(\bar{\xi}_s^f)^b$ are independent, and the processes do not jump at the same time, we obtain $[(\bar{\xi}_s^L)^{p-2}, \bar{\Theta}_s^g(\bar{\xi}_s^f)^b)]=0$ . Consider the differential term $d(\bar{\Theta}_s^g(\bar{\xi}_s^f)^b)$. Since $\bar{\xi}^f$ is a pure-jump process, we obtain
\begin{align*}
    \bar{\Theta}_s^g(\bar{\xi}_s^f)^b &= \int_0^s(\bar{\xi}_{t^-}^f)^bd\bar{\Theta}_t^g + \int_0^s\bar{\Theta}_{t^-}^gd(\bar{\xi}_{t}^f)^b+ [\bar{\Theta}_s^g,(\bar{\xi}_s^f)^b]  \\
    &= \int_0^s(\bar{\xi}_{t^-}^f)^bd\bar{\Theta}_t^g + \int_0^s\bar{\Theta}_{t^-}^gd(\bar{\xi}_{t}^f)^b+ \sum_{t\leq s} \Delta\bar{\Theta}_t^g \Delta(\bar{\xi}_t^f)^b \nonumber\\
    &= \int_0^s(\bar{\xi}_{t^-}^f)^bd\bar{\Theta}_t^g + \int_0^s\bar{\Theta}_{t^-}^gd(\bar{\xi}_{t}^f)^b+ \int_0^s \Delta\bar{\Theta}_t^g d(\bar{\xi}_s^f)^b \nonumber\\
    &= \int_0^s(\bar{\xi}_{t^-}^f)^bd\bar{\Theta}_t^g + \int_0^s\bar{\Theta}_{t}^gd(\bar{\xi}_{t}^f)^b
\end{align*}
which can be written in differential form as
\begin{equation*}
    d(\bar{\Theta}_s^g(\bar{\xi}_s^f)^b)= (\bar{\xi}_{s^-}^f)^bd\bar{\Theta}_s^g + \bar{\Theta}_s^gd(\bar{\xi}_{s}^f)^b
\end{equation*}
As a result, we can write $\mathbf{I_1}$ as\begin{align}
    \mathbf{I_1}&= \bar{\Theta}_t^g(\bar{\xi}_t^L)^{p-2}(\bar{\xi}_t^f)^b\int_0^t\sigma^2_2(\bar{\Theta}_s)ds 
    - \int_0^t\left(\int_0^s \sigma_2^2(\bar{\Theta}_{u^-})du\right)(\bar{\xi}_{s^-}^L)^{p-2}(\bar{\xi}_{s^-}^f)^bd\bar{\Theta}_s^g \label{2.}\\
    &- \int_0^t\left(\int_0^s \sigma_2^2(\bar{\Theta}_{u^-})du\right)(\bar{\xi}_{s^-}^L)^{p-2}\bar{\Theta}_s^gd(\bar{\xi}_{s}^f)^b
    - \int_0^t\left(\int_0^s \sigma_2^2(\bar{\Theta}_{u^-})du\right)\bar{\Theta}_{s^-}^g(\bar{\xi}_{s^-}^f)^bd(\bar{\xi}_s^L)^{p-2} \nonumber
    \end{align}
By similar arguments, $\mathbf{I_2}$ can be written as
\begin{align}
    \mathbf{I_2}&= \bar{\Theta}_t^{g-2}(\bar{\xi}_t^L)^{p}(\bar{\xi}_t^f)^b\int_0^t\sigma^2_2(\bar{\Theta}_s)ds 
    - \int_0^t\left(\int_0^s \sigma_2^2(\bar{\Theta}_{u^-})du\right)(\bar{\xi}_{s^-}^L)^{p}(\bar{\xi}_{s^-}^f)^bd\bar{\Theta}_s^{g-2}\label{3.}\\
    &- \int_0^t\left(\int_0^s \sigma_2^2(\bar{\Theta}_{u^-})du\right)(\bar{\xi}_{s^-}^L)^{p}\bar{\Theta}_s^{g-2}d(\bar{\xi}_{s}^f)^b
    - \int_0^t\left(\int_0^s \sigma_2^2(\bar{\Theta}_{u^-})du\right)\bar{\Theta}_{s^-}^{g-2}(\bar{\xi}_{s^-}^f)^bd(\bar{\xi}_s^L)^{p}\nonumber
    \end{align}
Next, we consider the term $\bar{\Theta}_s^g(\bar{\xi}_s^L)^p(\bar{\xi}_s^f)^b$ in $\mathbf{I_3}$, which can be written as\begin{align*}
    \bar{\Theta}_s^g(\bar{\xi}_s^L)^p(\bar{\xi}_s^f)^b&= (\bar{\Theta}_{s^-}+\Delta\bar{\Theta}_s)^g(\bar{\xi}_{s^-}^L+\Delta\bar{\xi}_s^L)^p(\bar{\xi}_{s^-}^f + \Delta\bar{\xi}_s^f)^b=: \mathbf{I_4}\cdot\mathbf{I_5}\cdot\mathbf{I_6}.
\end{align*}
Using the binomial expansions, we get
\begin{align*}
    \mathbf{I_4} &= \bar{\Theta}_{s^-}^g + \sum_{m_1=1}^g{g\choose m_1}\bar{\Theta}_{s^-}^{g-m_1}(\Delta\bar{\Theta}_s)^{m_1} \\
    \mathbf{I_5} &= (\bar{\xi}_{s^-}^L)^p + \sum_{m_2=1}^p{p\choose m_2}(\bar{\xi}_{s^-}^L)^{p-m_2}(\Delta\bar{\xi}_s^L)^{m_2} \\
    \mathbf{I_6} &= (\bar{\xi}_{s^-}^f)^b + \sum_{m_3=1}^b{b\choose m_3}(\bar{\xi}_{s^-}^f)^{b-m_3}(\Delta\bar{\xi}_s^f)^{m_3}.
\end{align*}
Then,\begin{align*}
     \mathbf{I_5}\cdot  \mathbf{I_6} &= (\bar{\xi}_{s^-}^L)^p(\bar{\xi}_{s^-}^f)^b +\sum_{m_2=1}^p{p\choose m_2}(\bar{\xi}_{s^-}^f)^b(\bar{\xi}_{s^-}^L)^{p-m_2}(\Delta\bar{\xi}_s^L)^{m_2} \\
     &+ \sum_{m_3=1}^b{b\choose m_3}(\bar{\xi}_{s^-}^L)^p(\bar{\xi}_{s^-}^f)^{b-m_3}(\Delta\bar{\xi}_s^f)^{m_3} \\
     &+ \sum_{m_2=1}^p{p\choose m_2}(\bar{\xi}_{s^-}^L)^{p-m_2}(\Delta\bar{\xi}_s^L)^{m_2}\sum_{m_3=1}^b{b\choose m_3}(\bar{\xi}_{s^-}^f)^{b-m_3}(\Delta\bar{\xi}_s^f)^{m_3}.
\end{align*}
Since $\bar{\xi}^L$ and $\bar{\xi}^f$ do not jump at the same time, the last term above vanishes. Then, we obtain\begin{align*}
    \mathbf{I_4}\cdot\mathbf{I_5}\cdot\mathbf{I_6} &= \bar{\Theta}_{s^-}^g(\bar{\xi}_{s^-}^L)^p(\bar{\xi}_{s^-}^f)^b +\sum_{m_2=1}^p{p\choose m_2}\bar{\Theta}_{s^-}^g(\bar{\xi}_{s^-}^f)^b(\bar{\xi}_{s^-}^L)^{p-m_2}(\Delta\bar{\xi}_s^L)^{m_2}\\
    &+ \sum_{m_3=1}^b{b\choose m_3}\bar{\Theta}_{s^-}^g(\bar{\xi}_{s^-}^L)^p(\bar{\xi}_{s^-}^f)^{b-m_3}(\Delta\bar{\xi}_s^f)^{m_3} \\
    &+ \sum_{m_1=1}^g{g\choose m_1}\bar{\Theta}_{s^-}^{g-m_1}(\bar{\xi}_{s^-}^L)^p(\bar{\xi}_{s^-}^f)^b(\Delta\bar{\Theta}_s)^{m_1} \\
    &+ \sum_{m_1=1}^g{g\choose m_1}\bar{\Theta}_{s^-}^{g-m_1}(\Delta\bar{\Theta}_s)^{m_1}\sum_{m_2=1}^p{p\choose m_2}(\bar{\xi}_{s^-}^f)^b(\bar{\xi}_{s^-}^L)^{p-m_2}(\Delta\bar{\xi}_s^L)^{m_2} \\
    &+\sum_{m_1=1}^g{g\choose m_1}\bar{\Theta}_{s^-}^{g-m_1}(\Delta\bar{\Theta}_s)^{m_1}\sum_{m_3=1}^b{b\choose m_3}(\bar{\xi}_{s^-}^L)^p(\bar{\xi}_{s^-}^f)^{b-m_3}(\Delta\bar{\xi}_s^f)^{m_3}
\end{align*}
The fifth term vanishes as $\bar{\Theta}$ and $\bar{\xi}^L$ do not jump at the same time. For the second, third, and fourth terms, we can write the terms corresponding to $m_1=1$, $m_2=1$, and $m_3=1$ separated from the sums to obtain that \begin{align*}
    \mathbf{I_4}\cdot\mathbf{I_5}\cdot\mathbf{I_6} &= \bar{\Theta}_{s^-}^g(\bar{\xi}_{s^-}^L)^p(\bar{\xi}_{s^-}^f)^b +\, p\bar{\Theta}_{s^-}^g(\bar{\xi}_{s^-}^f)^b(\bar{\xi}_{s^-}^L)^{p-1}\Delta\bar{\xi}_s^L\\ 
    &+\, b\bar{\Theta}_{s^-}^g(\bar{\xi}_{s^-}^L)^p(\bar{\xi}_{s^-}^f)^{b-1}\Delta\bar{\xi}_s^f+g\bar{\Theta}_{s^-}^{g-1}(\bar{\xi}_{s^-}^L)^p(\bar{\xi}_{s^-}^f)^b\Delta\bar{\Theta}_s\\
   &+\, \sum_{m_2=2}^p{p\choose m_2}\bar{\Theta}_{s^-}^g(\bar{\xi}_{s^-}^f)^b(\bar{\xi}_{s^-}^L)^{p-m_2}(\Delta\bar{\xi}_s^L)^{m_2}\\
    &+\, \sum_{m_3=2}^b{b\choose m_3}\bar{\Theta}_{s^-}^g(\bar{\xi}_{s^-}^L)^p(\bar{\xi}_{s^-}^f)^{b-m_3}(\Delta\bar{\xi}_s^f)^{m_3} \\
    &+ \sum_{m_1=2}^g{g\choose m_1}\bar{\Theta}_{s^-}^{g-m_1}(\bar{\xi}_{s^-}^L)^p(\bar{\xi}_{s^-}^f)^b(\Delta\bar{\Theta}_s)^{m_1} \\
    &+\,\sum_{m_1=1}^g{g\choose m_1}\bar{\Theta}_{s^-}^{g-m_1}(\Delta\bar{\Theta}_s)^{m_1}\sum_{m_3=1}^b{b\choose m_3}(\bar{\xi}_{s^-}^L)^p(\bar{\xi}_{s^-}^f)^{b-m_3}(\Delta\bar{\xi}_s^f)^{m_3}
\end{align*}
Hence we have\begin{align*}
    \mathbf{I_3} &= \sum_{0\leq s \leq t}[ \sum_{m_2=2}^p{p\choose m_2}\bar{\Theta}_{s^-}^g(\bar{\xi}_{s^-}^f)^b(\bar{\xi}_{s^-}^L)^{p-m_2}(\Delta\bar{\xi}_s^L)^{m_2}\\
    &+ \sum_{m_3=2}^b{b\choose m_3}\bar{\Theta}_{s^-}^g(\bar{\xi}_{s^-}^L)^p(\bar{\xi}_{s^-}^f)^{b-m_3}(\Delta\bar{\xi}_s^f)^{m_3} \\
    &+ \sum_{m_1=2}^g{g\choose m_1}\bar{\Theta}_{s^-}^{g-m_1}(\bar{\xi}_{s^-}^L)^p(\bar{\xi}_{s^-}^f)^b(\Delta\bar{\Theta}_s)^{m_1} \\
    &+\sum_{m_1=1}^g{g\choose m_1}\bar{\Theta}_{s^-}^{g-m_1}(\Delta\bar{\Theta}_s)^{m_1}\sum_{m_3=1}^b{b\choose m_3}(\bar{\xi}_{s^-}^L)^p(\bar{\xi}_{s^-}^f)^{b-m_3}(\Delta\bar{\xi}_s^f)^{m_3}]
\end{align*}
We can write $\mathbf{I_3}$ as a sum of some stochastic integrals.\begin{align*}
     \mathbf{I_3} &= \sum_{m_2=2}^p\int_0^t{p\choose m_2}\bar{\Theta}_{s^-}^g(\bar{\xi}_{s^-}^f)^b(\bar{\xi}_{s^-}^L)^{p-m_2}d\xi_s^{(m_2)}\\
    &+ \sum_{m_3=2}^b\int_0^t{b\choose m_3}\bar{\Theta}_{s^-}^g(\bar{\xi}_{s^-}^L)^p(\bar{\xi}_{s^-}^f)^{b-m_3}d\xi_s^{[m_3]} \\
    &+ \sum_{m_1=2}^g\int_0^t{g\choose m_1}\bar{\Theta}_{s^-}^{g-m_1}(\bar{\xi}_{s^-}^L)^p(\bar{\xi}_{s^-}^f)^bd\Theta_s^{(m_1)} \\
    &+\sum_{m_1=1}^g\int_0^t{g\choose m_1}\bar{\Theta}_{s^-}^{g-m_1}d\Theta_s^{(m_1)}\sum_{m_3=1}^b\int_0^t{b\choose m_3}(\bar{\xi}_{s^-}^L)^p(\bar{\xi}_{s^-}^f)^{b-m_3}d\xi_s^{[m_3]}
\end{align*}
We should write this expression in terms of the compensated processes $\bar{\Theta}^{(k)},\bar{\xi}^{(l)},\bar{\xi}^{[m]}$. 
\begin{align}
     \mathbf{I_3} &= \sum_{m_2=2}^p\int_0^t{p\choose m_2}\bar{\Theta}_{s^-}^g(\bar{\xi}_{s^-}^f)^b(\bar{\xi}_{s^-}^L)^{p-m_2}d\bar{\xi}_s^{(m_2)}\nonumber\\
    &+ \sum_{m_3=2}^b\int_0^t{b\choose m_3}\bar{\Theta}_{s^-}^g(\bar{\xi}_{s^-}^L)^p(\bar{\xi}_{s^-}^f)^{b-m_3}d\bar{\xi}_s^{[m_3]}\nonumber \\
    &+ \sum_{m_1=2}^g\int_0^t{g\choose m_1}{\bar{\Theta}}_{s^-}^{g-m_1}(\bar{\xi}_{s^-}^L)^p(\bar{\xi}_{s^-}^f)^bd\bar{\Theta}_s^{(m_1)} \nonumber\\
    &+\sum_{m_1=1}^g\int_0^t{g\choose m_1}\bar{\Theta}_{s^-}^{g-m_1}d\bar{\Theta}_s^{(m_1)}\sum_{m_3=1}^b\int_0^t{b\choose m_3}(\bar{\xi}_{s^-}^L)^p(\bar{\xi}_{s^-}^f)^{b-m_3}d\bar{\xi}_s^{[m_3]} + \mathbf{I_7} \label{4.}
\end{align}
where 
\begin{align*}
     \mathbf{I_7}& = \sum_{m_2=2}^p\int_0^t\int_\mathbb{R}{p\choose m_2}\bar{\Theta}_{s^-}^g(\bar{\xi}_{s^-}^f)^b(\bar{\xi}_{s^-}^L)^{p-m_2}x^{m_2}\nu_2(dx)ds\\
    &+ \sum_{m_3=2}^b\int_0^t\int_\mathbb{R}\int_\mathbb{R}{b\choose m_3}\bar{\Theta}_{s^-}^g(\bar{\xi}_{s^-}^L)^p(\bar{\xi}_{s^-}^f)^{b-m_3}y^{m_3}\nu_1(dx)\mathbb{P}(U(x)\in dy)ds\\
    &+ \sum_{m_1=2}^g\int_0^t\int_\mathbb{R}{g\choose m_1}{\bar{\Theta}}_{s^-}^{g-m_1}(\bar{\xi}_{s^-}^L)^p(\bar{\xi}_{s^-}^f)^b\nu_1(dx)ds \\
    &+\sum_{m_1=1}^g\int_0^t\int_\mathbb{R}{g\choose m_1}\bar{\Theta}_{s^-}^{g-m_1}\nu_1(dx)ds\sum_{m_3=1}^b\int_0^t{b\choose m_3}(\bar{\xi}_{s^-}^L)^p(\bar{\xi}_{s^-}^f)^{b-m_3}d\bar{\xi}_s^{[m_3]}\\
    &+ \sum_{m_1=1}^g\int_0^t{g\choose m_1}\bar{\Theta}_{s^-}^{g-m_1}d\bar{\Theta}_s^{(m_1)}\\
    &\cdot\sum_{m_3=1}^b\int_0^t\int_\mathbb{R}\int_\mathbb{R}{b\choose m_3}(\bar{\xi}_{s^-}^L)^p(\bar{\xi}_{s^-}^f)^{b-m_3}y^{m_3}\nu_1(dx)\mathbb{P}(U(x) \in dy)ds\\
    &+ \sum_{m_1=1}^g\int_0^t\int_\mathbb{R}{g\choose m_1}\bar{\Theta}_{s^-}^{g-m_1}\nu_1(dx)ds\\
    &\cdot \sum_{m_3=1}^b\int_0^t\int_\mathbb{R}\int_\mathbb{R}{b\choose m_3}(\bar{\xi}_{s^-}^L)^p(\bar{\xi}_{s^-}^f)^{b-m_3}y^{m_3}\nu_1(dx)\mathbb{P}(U(x) \in dy)ds
\end{align*}
Note that $\mathbf{I_7}$ consists of terms of the form\begin{equation*}
    \int_0^t \bar{\Theta}_{s^-}^{g-m_1}(\bar{\xi}_{s^-}^L)^{p-m_2}(\bar{\xi}_{s^-}^f)^{b-m_3}\kappa(s^-)ds
\end{equation*} 
for various predictable processes $\kappa(s^-)$. Then, as in $\mathbf{I_1}$, we can write $\mathbf{I_7}$ as a sum of stochastic integrals with respect to processes $\bar{\Theta}^k$, $(\bar{\xi}^L)^l$, and $(\bar{\xi}^f)^m$ with lower powers. Indeed, by integration by parts, we obtain \begin{align*}
    \int_0^t \bar{\Theta}_{s^-}^{g-m_1}(\bar{\xi}_{s^-}^L)^{p-m_2}(\bar{\xi}_{s^-}^f)^{b-m_3}\kappa(s^-)ds &= \bar{\Theta}_t^{g-m_1}(\bar{\xi}_t^L)^{p-m_2}(\bar{\xi}_t^f)^{b-m_3}\int_0^t\kappa(s^-)ds \\
    &- \int_0^t\left(\int_0^s \kappa(u^-)du\right)(\bar{\xi}_{s^-}^L)^{p-m_2}(\bar{\xi}_{s^-}^f)^{b-m_3}d\bar{\Theta}_s^{g-m_1}\\
    &- \int_0^t\left(\int_0^s \kappa(u^-)du\right)(\bar{\xi}_{s^-}^L)^{p-m_2}\bar{\Theta}_s^{g-m_1}d(\bar{\xi}_{s}^f)^{b-m_3}\\
    &- \int_0^t\left(\int_0^s \kappa(u^-)du\right)\bar{\Theta}_{s^-}^{g-m_1}(\bar{\xi}_{s^-}^f)^{b-m_3}d(\bar{\xi}_s^L)^{p-m_2}
\end{align*}
Writing $\bar{\Theta}=\bar{\Theta}^{(1)}$, $\bar{\xi}^L=\bar{\xi}^{(1)}$, $\bar{\xi}^f=\bar{\xi}^{[1]}$, and  combining Equations \eqref{1.},\eqref{2.},\eqref{3.}, and \eqref{4.}, we get\begin{align*}
    \bar{\Theta}_t^g &(\bar{\xi}_t^L)^p (\bar{\xi}_t^f)^b = 
    \bar{\Theta}_0^g(\bar{\xi}_0^L)^p(\bar{\xi}_0^f)^b + \int_{0^+}^t g\bar{\Theta}_{s^-}^{g-1}(\bar{\xi}_{s^-
    }^L)^p(\bar{\xi}_{s^-}^f)^b d\bar{\Theta}_s^{(1)} + \int_{0^+}^tp\bar{\Theta}_{s^-}^g (\bar{\xi}_{s^-}^L)^{p-1}
    (\bar{\xi}_{s^-}^f)^bd\bar{\xi}_s^{(1)}   \\
    &+ \int_{0^+}^t b\bar{\Theta}_{s^-}^g(\bar{\xi}_{s^-}^L)^p(\bar{\xi}_{s^-}^f)^{b-1}d\bar{\xi}_{s}^{[1]} + \frac{1}{2}p(p-1)\left[ \bar{\Theta}_t^g(\bar{\xi}_t^L)^{p-2}(\bar{\xi}_t^f)^b\int_0^t\sigma^2_2(\bar{\Theta}_s)ds\right. \\
    &- \int_0^t(\int_0^s \sigma_2^2(\bar{\Theta}_{u^-})du)(\bar{\xi}_{s^-}^L)^{p-2}(\bar{\xi}_{s^-}^f)^bd\bar{\Theta}_s^g- \int_0^t(\int_0^s \sigma_2^2(\bar{\Theta}_{u^-})du)(\bar{\xi}_{s^-}^L)^{p-2}\bar{\Theta}_s^gd(\bar{\xi}_{s}^f)^b\\
    &\left. - \int_0^t(\int_0^t \sigma_2^2(\bar{\Theta}_{u^-})du)\bar{\Theta}_{s^-}^g(\bar{\xi}_{s^-}^f)^bd(\bar{\xi}_s^L)^{p-2} \right] +\frac{1}{2}g(g-1)\left[\bar{\Theta}_t^{g-2}(\bar{\xi}_t^L)^{p}(\bar{\xi}_t^f)^b\int_0^t\sigma^2_2(\bar{\Theta}_s)ds \right. \\
    &- \int_0^t(\int_0^s \sigma_2^2(\bar{\Theta}_{u^-})du)(\bar{\xi}_{s^-}^L)^{p}(\bar{\xi}_{s^-}^f)^bd\bar{\Theta}_s^{g-2}- \int_0^t(\int_0^s \sigma_2^2(\bar{\Theta}_{u^-})du)(\bar{\xi}_{s^-}^L)^{p}\bar{\Theta}_s^{g-2}d(\bar{\xi}_{s}^f)^b\\
    &\left. - \int_0^t(\int_0^s \sigma_2^2(\bar{\Theta}_{u^-})du)\bar{\Theta}_{s^-}^{g-2}(\bar{\xi}_{s^-}^f)^bd(\bar{\xi}_s^L)^{p}\right] + \sum_{m_2=2}^p\int_0^t{p\choose m_2}\bar{\Theta}_{s^-}^g(\bar{\xi}_{s^-}^f)^b(\bar{\xi}_{s^-}^L)^{p-m_2}d\bar{\xi}_s^{(m_2)}\nonumber\\
    &+ \sum_{m_3=2}^b\int_0^t{b\choose m_3}\bar{\Theta}_{s^-}^g(\bar{\xi}_{s^-}^L)^p(\bar{\xi}_{s^-}^f)^{b-m_3}d\bar{\xi}_s^{[m_3]} + \sum_{m_1=2}^g\int_0^t{g\choose m_1}{\bar{\Theta}}_{s^-}^{g-m_1}(\bar{\xi}_{s^-}^L)^p(\bar{\xi}_{s^-}^f)^bd\bar{\Theta}_s^{(m_1)} \nonumber\\
    &+\sum_{m_1=1}^g\int_0^t{g\choose m_1}\bar{\Theta}_{s^-}^{g-m_1}d\bar{\Theta}_s^{(m_1)}\sum_{m_3=1}^b\int_0^t{b\choose m_3}(\bar{\xi}_{s^-}^L)^p(\bar{\xi}_{s^-}^f)^{b-m_3}d\bar{\xi}_s^{[m_3]} + \mathbf{I_7}
\end{align*}
As a result, we write $\Theta^g (\xi^L)^p (\xi^f)^b $ as a sum of stochastic integrals of powers of $\bar{\Theta}$, $\bar{\xi}^f$, and $\bar{\xi}^L$ strictly less than $g+p+b$ with respect to processes $\bar{\Theta}^{(k)}$, $\bar{\xi}^{(l)}$, and $\bar{\xi}^{[m]}$. Hence by induction, we obtain the desired result. \hfill $\Box$\\

Now, we will show that the polynomials of $\Theta,\xi^L$ and $\xi^f$ are in dense in $L^2(\Omega,\mathcal{F})$ in the next proposition. 
\begin{Prop}\label{Base}
    For fixed $t\geq 0$ and $0\leq s_1 \leq \ldots \leq s_m\leq t$, let\begin{align*}
        \mathbf{\Theta} &= (\bar{\Theta}_{s_1},\ldots, \bar{\Theta}_{s_m})\\
        \mathbf{\xi^L} &= (\bar{\xi}_{s_1}^L,\ldots , \bar{\xi}_{s_m}^L) \\
        \mathbf{\xi^f} &= (\bar{\xi}_{s_1}^f,\ldots,\bar{\xi}_{s_m}^f)
    \end{align*}
and let $Z\in L^2(\Omega,\mathcal{F}_t)$. Then, for any $\epsilon>0$, there exist an integer $m$ and a random variable $Z_\epsilon \in L^2(\Omega,\sigma(\mathbf{\Theta},\mathbf{\xi^L},\mathbf{\xi^f}))$ such that\begin{equation*}
    \mathbb{E}[(Z-Z_\epsilon)^2]< \epsilon
\end{equation*}
Furthermore, the following set is a total family in $L^2(\Omega,\mathcal{F})$\begin{align*}
    \mathcal{P} =& \{ \bar{\Theta}_{t_1}^{g_1}\cdot \ldots \cdot \bar{\Theta}_{t_m}^{g_m}\cdot (\bar{\xi}_{t_1}^L)^{p_1}\cdot \ldots \cdot (\bar{\xi}_{t_m}^L)^{p_m}\cdot (\bar{\xi}_{t_1}^f)^{b_1}\cdot \ldots \cdot (\bar{\xi}_{t_m}^f)^{b_m} \, \\&: \, 0\leq t_1<\ldots <t_m; g_1,\ldots,g_m,p_1,\ldots,p_m,b_1,\ldots,b_m\geq 0\}.
\end{align*}
\end{Prop}
\textbf{Proof:} The result follows from the general considerations as in the proof of \cite[Lem.4.3.1]{MR2001996}. Indeed, let $\{ t_i\}_{i=1}^\infty$ be a dense subset of $[0,t]$, and $\mathcal{H}_n$ be the $\sigma$-algebra generated by $\bar{\Theta}_{t_1},\ldots,\bar{\Theta}_{t_n},\bar{\xi}^L_{t_1},\ldots,\bar{\xi}^L_{t_n},\bar{\xi}^f_{t_1},\ldots,\bar{\xi}^f_{t_n}$. Then, $\mathcal{H}_n\subset \mathcal{H}_{n+1}$ for each $n\geq 1$, and $\mathcal{F}_t$ is the smallest $\sigma$-algebra containing all the $\mathcal{H}_n$'s. Let $Z\in L^2(\Omega, \mathcal{F}_t)$. By martingale convergence theorem \cite[Cor.C.9]{MR2001996}, we have \begin{equation*}
    Z = \mathbb{E}[Z|\mathcal{F}_t]= \lim_{n\to\infty}\mathbb{E}[Z|\mathcal{H}_n]
\end{equation*}
where the limit is in $L^2(\Omega,\mathcal{F}_t)$ and pointwise. It follows that there exists $m\in \mathbb{N}$ such that $\mathbb{E}[(Z-\mathbb{E}[Z|\mathcal{H}_m])^2]<\epsilon$ ,and we let $Z_\epsilon =\mathbb{E}[Z|\mathcal{H}_m]\in L^2(\Omega,\sigma(\mathbf{\Theta,\xi^L,\xi^f}))$. For the second part, assume that there exists $Z\in L^2(\Omega,\mathcal{F})$ which is orthogonal to $\mathcal{P}$. By the first part of proof, there exist $Z_\epsilon\in L^2(\Omega,\sigma(\mathbf{\Theta},\mathbf{\xi^L},\mathbf{\xi^f}))$ such that $\mathbb{E}[(Z-Z_\epsilon)^2]< \epsilon$. Note that polynomials are dense in $L^2(\mathbb{R},d\varphi)$ by the assumptions given at the end of Section \ref{Preliminaries}. Therefore, we can approximate $Z_\epsilon$ by polynomials from $\mathcal{P}$ and we get $\mathbb{E}[Z_\epsilon Z]=0$. Then, we have\begin{equation*}
    \mathbb{E}[Z^2] = \mathbb{E}[Z^2-Z_\epsilon Z] = \mathbb{E}[Z(Z-Z_\epsilon)]\leq \sqrt{\mathbb{E}[Z^2]\mathbb{E}[(Z-Z_\epsilon)^2]}\leq \sqrt{\epsilon \mathbb{E}[Z^2]}
\end{equation*}
by Schwarz inequality. As $\epsilon \to 0$, we obtain $Z=0$ a.s. \hfill $\Box$

Now, we are ready to give the chaotic representation property in the next theorem. It states that any square-integrable random variable can be represented as a stochastic integral with respect to pairwise orthogonal martingales.
\begin{Thm}\label{Chaotic Representation Theorem}
    Any square-integrable $(\mathcal{F}_t)$-measurable random variable F can be represented as follows:\begin{align*}
        F &= \mathbb{E}[F] + \sum_{s=1}^\infty \sum_{\tau=1}^\infty \sum_{\zeta=1}^\infty \sum_{\nu_1,\ldots,\nu_\zeta\geq 1} \sum_{\mu_1,\ldots,\mu_\tau\geq 0 } \sum_{\iota_1,\ldots,\iota_s\geq 1}
         \int_0^t\int_0^{t_{1^-}}\ldots \int_0^{t_{\zeta+s+\tau-1}^-}\\ &f_{\nu_1,\ldots,\nu_\zeta,\mu_1,\ldots,\mu_\tau,\iota_1,\ldots,\iota_s}(t,t_1,t_2,\ldots,t_{\zeta+s+\tau})
         dH_{t_\zeta+s+\tau}^{(\nu_\zeta)}\ldots dH_{t_{1+s+\tau}}^{(\nu_1)}dG_{t_{s+\tau}}^{(\iota_s,\mu_\tau)}\ldots dG_{t_1}^{(\iota_1,\mu_0)}
    \end{align*}
    where $f_{\nu_1,\ldots,\nu_\zeta,\mu_1,\ldots,\mu_\tau,\iota_1,\ldots,\iota_s}(t,t_1,t_2,\ldots,t_{\zeta+s+\tau})$ are some random field such that the stochastic integral is well-defined on $L^2(\Omega,\mathcal{F})$, $H^{(k)}$ and $G^{(l,m)}$ for $k\geq 1,l\geq 1,m\geq 0$ are orthogonal martingales given in \eqref{H} and \eqref{G}, respectively. 
\end{Thm}
\begin{Rem}\label{Remark}
   The series in Theorem \ref{Chaotic Representation Theorem} is the $L^2(\Omega,\mathcal{F})$-limit of
    \begin{equation*}
        \sum_{s=1}^{n} \sum_{\tau=1}^{q} \sum_{\zeta=1}^{r} \sum_{\nu_1,\ldots,\nu_\zeta\geq 1} \sum_{\mu_1,\ldots,\mu_\tau\geq 0 } \sum_{\iota_1,\ldots,\iota_s\geq 1}
         \int_0^t\int_0^{t_{1^-}}\ldots \int_0^{t_{\zeta+s+\tau-1}^-}\ldots
    \end{equation*}
    as $n,q,r\to \infty$\cite[Thm.2,Remark 3]{MR3854533}.
\end{Rem}
\textbf{Proof:} Using Proposition \ref{Representation of Powers}, and linear transformations from $\bar{\xi}^{(k)}$ to $H^{(k)}$, from $\bar{\Theta}^{(l)}+\bar{\xi}^{[m]}$ to $G^{(l,m+1)}$, we can write
    \begin{align}
        \bar{\Theta}_t^g& (\bar{\xi}_t^L)^p (\bar{\xi}_t^f)^b = f^{(g+p+b)}(t) + \sum_{s=1}^p \sum_{\tau=1}^g \sum_{\zeta=1}^b \sum_{\nu_1,\ldots,\nu_\zeta\geq 1} \sum_{\mu_1,\ldots,\mu_\tau\geq 0 } \sum_{\iota_1,\ldots,\iota_s\geq 1}
         \int_0^t\int_0^{t_{1^-}}\ldots \int_0^{t_{\zeta+s+\tau-1}^-}\nonumber\\ &f_{\nu_1,\ldots,\nu_\zeta,\mu_1,\ldots,\mu_\tau,\iota_1,\ldots,\iota_s}^{(g+p+b)}(t,t_1,t_2,\ldots,t_{\zeta+s+\tau}) dH_{t_\zeta+s+\tau}^{(\nu_\zeta)}\ldots dH_{t_{1+s+\tau}}^{(\nu_1)}dG_{t_{s+\tau}}^{(\iota_s,\mu_\tau)}\ldots dG_{t_1}^{(\iota_1,\mu_1)}\label{Representation}
    \end{align}
where $f^{(g+p+b)}(t)$ and $f_{(\nu_1,\ldots , \nu_\zeta,\mu_1,\ldots , \mu_\tau, \iota_1, \ldots, \iota_s)}^{(g+p+b)}(t,t_1,\ldots, t_{\zeta+s+\tau})$ are sum of product of predictable processes observed at $t,t_1,\ldots, t_{\zeta+s+\tau -1}$ defined on $L^2(\Omega,\mathcal{F})$. Since $\mathcal{P}$ is a total family by Proposition \ref{Base}, we can represent any square integrable $(\mathcal{F}_t)$-measurable random variable $F$ in terms of $\bar{\Theta}_t^g (\bar{\xi}_t^L)^p (\bar{\xi}_t^f)^b$ for $p,g,b\geq 0$. As in the proof of \cite[Thm.1]{MR1787127}, this representation can be written in the form of \eqref{Representation}. \hfill $\Box$\\

Finally, we give the predictable representation property for a MAP in the next theorem. It states that any square-integrable $(\mathcal{F}_t)$-martingale $M$ can be represented as a sum of stochastic integrals of predictable processes with respect to the processes $\bar{\xi}$, $\bar{\xi}^{(k)}$, $\bar{\xi}^{[l]}$ and $\bar{\Theta}^{(m)}$ for $k \geq 2$, and $l,m\geq 1$ where $\bar{\xi}$ is the compensated version of $\xi$. The notation is the same as in Remark \ref{Remark}.
\begin{Thm}\label{Main Predictable Representation Theorem}
    Any square-integrable $(\mathcal{F}_t)$-martingale $M$ can be represented as\begin{equation*}
        M_t=M_0 +\int_0^t h_\xi(s)d\bar{\xi}_s + \sum_{k=2}^\infty \int_0^t h_L^{(k)}(s)(s)d\bar{\xi}_s^{(k)}+ \sum_{l=1}^\infty \int_0^t h_f^{(l)}(s)d\bar{\xi}_s^{[l]} 
        +\sum_{m=1}^\infty \int_0^t h_\Theta^{(m)}(s)d\bar{\Theta}_s^{(m)}
    \end{equation*}
where $h_\xi,h_L^{(k)},h_f^{(l)}$, and $h_\Theta^{(m)}$ are predictable processes.
\end{Thm}
\textbf{Proof:} As in Theorem \ref{Chaotic Representation Theorem}, we can switch back from $H^{(k)}$ to $\bar{\xi}^{(k)}$, and from $G^{(l,m)}$ to $\bar{\Theta}^{(l)}$ and $\bar{\xi}^{[m-1]}$ by linear transformations to obtain that\begin{align*}
       F(t)-\mathbb{E}[F(t)] &= \sum_{\zeta=1}^p \sum_{s=1}^g\sum_{\tau=1}^b\sum_{(\nu_1,\ldots,\nu_\zeta)\in \{1,\ldots,p\}^\zeta} \sum_{(\mu_1,\ldots,\mu_\tau)\in \{1,\ldots,b\}^\tau}\sum_{(\iota_1,\ldots,\iota_s)\in \{1,\ldots,g\}^s}\\
        & \int_0^t \int_0^{t_1^-} \ldots \int_0^{t_{\zeta + s+\tau -1}^-} f_{(\nu_1,\ldots , \nu_\zeta,\mu_1,\ldots , \mu_\tau, \iota_1, \ldots, \iota_s)}^{(g+p+b)}(t,t_1,\ldots, t_{\zeta+s+\tau})\\
        & \quad \quad \quad \quad \quad d\bar{\xi}_{t_{\zeta+s+\tau}}^{(\nu_\zeta)} \ldots d\bar{\xi}_{t_{1+s+\tau}}^{(\nu_1)}d\bar{\Theta}_{t_{s+\tau}}^{(\iota_s)}\ldots\bar{\Theta}_{t_{\tau+1}}^{(\iota_1)}d\bar{\xi}_{t_{\tau}}^{[\mu_\tau]}\ldots d\bar{\xi}_{t_1}^{[\mu_1]}
    \end{align*}
As in the proof of \cite[Thm.4]{MR3854533} and \cite[Thm.1]{MR1787127}, we get\begin{equation*}
        F_t-\mathbb{E}[F(t)]= \sum_{k=1}^\infty \int_0^t h_L^{(k)}(s)d\bar{\xi}_s^{(k)}+ \sum_{l=1}^\infty \int_0^t \tilde{h}_f^{(l)}(s)d\bar{\xi}_s^{[l]} +\sum_{m=1}^\infty \int_0^t h_\Theta^{(m)}(s)d\bar{\Theta}_s^{(m)}
    \end{equation*}
where $h_\xi,h_L^{(k)},\tilde{h}_f^{(l)}$, and $h_\Theta^{(m)}$ are predictable processes. By definition of $\bar{\xi}$, we have\begin{equation*}
    \int_0^t h_L^{(1)}(s)d\bar{\xi}_s^{(1)}= \int_0^1 h^{(1)}_L(s)d\bar{\xi}_s-\int_0^1{h}^{(1)}_L(s)d\bar{\xi}^{[1]}
\end{equation*}
    Then, by defining $h_\xi(s)=h_L^{(1)}(s)$, $h_f^{(1)}=\tilde{h}_f^{(1)}-h_L^{(1)}$, and $h_f^{(l)}=\tilde{h}_f^{(l)}$ for $l\geq 2$, we obtain the desired representation. \hfill $\Box$
    \\  
\textbf{Acknowledgments}\\
First author is supported by the Scientific and Technological Research Council of Turkey (TÜBİTAK) through BIDEB 2211.
\bibliographystyle{apalike}
\bibliography{references}

\end{document}